\documentclass[12pt,twoside]{article}
\usepackage[english]{babel}
\usepackage{umj_3_1}           
\usepackage{graphicx, srcltx}   
\usepackage{cite}               
\usepackage{url}                                

\usepackage{amsmath}            
\usepackage{amsfonts,amssymb}   
\usepackage{srcltx}             

\newcommand{\Leftclt}{Victor Nijimbere}   
\newcommand{\Rightclt}{Some non-elementary integrals of sine, cosine and exponential
integrals type} 
\usepackage{authblk}

\begin{document}

\date{}

\Mainclt            


\begin{Titul}
{\large \bf EVALUATION OF SOME NON-ELEMENTARY INTEGRALS INVOLVING SINE, COSINE, EXPONENTIAL AND LOGARITHMIC INTEGRALS: PART II}\\ \vspace{5mm}
{{\bf Victor Nijimbere} \\ [2ex] {\small School of Mathematics and Statistics, Carleton University, \hspace{7cm}Ottawa, Ontario, Canada,\hspace{9cm}victornijimbere@gmail.com}\\[3ex]} 
\end{Titul}

\begin{Anot}
{\bf Abstract:} The non-elementary integrals $\mbox{Si}_{\beta,\alpha}=\int [\sin{(\lambda x^\beta)}/(\lambda x^\alpha)] dx,\beta\ge1,\alpha>\beta+1$ and $\mbox{Ci}_{\beta,\alpha}=\int [\cos{(\lambda x^\beta)}/(\lambda x^\alpha)] dx, \beta\ge1, \alpha>2\beta+1$, where $\{\beta,\alpha\}\in\mathbb{R}$, are evaluated in terms of the hypergeometric function  $_{2}F_3$. On the other hand, the exponential integral $\mbox{Ei}_{\beta,\alpha}=\int (e^{\lambda x^\beta}/x^\alpha) dx, \beta\ge1, \alpha>\beta+1$ is expressed in terms of $_{2}F_2$. The method used to evaluate these integrals consists of expanding the integrand  as a Taylor series and integrating the series term by term. 


\end{Anot}

\section{Introduction}\label{sec:1}
\setcounter{equation}{0}

 Let us first give the definition of the non-elementary integral. This definition is also given in Part I \cite{NV1}, we repeat it here for reference.

\begin{defen} An elementary function is a function of one variable constructed using that variable and constants, and by performing a finite number of repeated algebraic operations involving exponentials and logarithms.
An indefinite integral which can be expressed in terms of elementary functions is an elementary integral. And if, on the other hand, it cannot be evaluated in terms of elementary functions, then it is non-elementary \cite{MZ,R}.
\label{defen1}
\end{defen}

 The cases consisting of the non-elementary integrals $\mbox{Si}_{\beta,\alpha}=\int [\sin{(\lambda x^\beta)}/(\lambda x^\alpha)] dx, \beta\ge1,\alpha\le\beta+1$ and $\mbox{Ci}_{\beta,\alpha}=\int [\cos{(\lambda x^\beta)}/(\lambda x^\alpha)] dx$, $\beta\ge1, \alpha\le2\beta+1$, where $\{\beta,\alpha\}\in\mathbb{R}$, were considered and evaluated in terms of the hypergeometric functions $_{1}F_2$ and $_{2}F_3$ in Part I \cite{NV1}, and their asymptotic expressions for $|x|\gg1$ were derived too in Part I \cite{NV1}.  
The exponential integral $\mbox{Ei}_{\beta,\alpha}=\int (e^{\lambda x^\beta}/x^\alpha) dx$ where $\beta\ge1$ and $\alpha\le\beta+1$ was expressed in terms of $_{2}F_2$, and its asymptotic expression for $|x|\gg1$ was derived as well in Part I \cite{NV1}. 

 Here, we investigate other cases which were not treated neither in Part I \cite {NV1} nor elsewhere. We evaluate $\mbox{Si}_{\beta,\alpha}=\int [\sin{(\lambda x^\beta)}/(\lambda x^\alpha)] dx, \beta\ge1,\alpha>\beta+1$ and $\mbox{Ci}_{\beta,\alpha}=\int [\cos{(\lambda x^\beta)}/(\lambda x^\alpha)] dx$, $\beta\ge1, \alpha>2\beta+1$ and $\mbox{Ei}_{\beta,\alpha}=\int (e^{\lambda x^\beta}/x^\alpha)dx, \beta\ge1, \alpha>\beta+1$.
In order to take into account all possibilities, we write these integrals  as $\mbox{Si}_{\beta,\beta+\alpha}=\int [\sin{(\lambda x^\beta)}/(\lambda x^{\beta+\alpha})] dx, \beta\ge1, \alpha>1$, $\mbox{Ci}_{\beta,2\beta+\alpha}=\int [\cos{(\lambda x^\beta)}/(\lambda x^{2\beta+\alpha})] dx$, $\beta\ge1, \alpha>1$, and $\mbox{Ei}_{\beta,\beta+\alpha}=\int (e^{\lambda x^\beta}/x^{\beta+\alpha}) dx, \beta\ge1, \alpha>1$ where $\{\beta,\alpha\}\in\mathbb{R}$. On one hand, $\mbox{Si}_{\beta,\beta+\alpha}$ and $\mbox{Ci}_{\beta,2\beta+\alpha}$ are expressed in terms of  the hypergeometric function $_{2}F_3$, while on another hand, $\mbox{Ei}_{\beta,\beta+\alpha}$ is expressed in terms of the hypergeometric function $_{2}F_2$. 


These integrals involving a power function $ x^\beta$ in the argument of the numerator are the generalizations of the exponential, sine and cosine integrals in \cite{N} (see sections 8.19 and 8.21 respectively), which have applications in different fields in science, applied sciences and engineering including physics, nuclear technology, mathematics, probability, statistics, and so on. For instance, the generalized exponential integral $E_{1,1+\alpha}$ is used in fluidodynamics and transport theory, where it is applied to the solution of Milne's integral equations \cite{C}, there are also used in modeling radiative transfer processes in the atmosphere and in nuclear reactors \cite{S}, etc. Exponential asymptotics involving generalized exponential integrals are used in probability theory, see for example \cite{Ch}. On the hand, generalized sine and cosine integrals are frequently utilized in Fourier analysis and related domains \cite{Ra}. Therefore, we are justified to further generalize these functions and their connections to hypergeometric functions.


Before we proceed to the main objectives of this paper consisting of evaluating the above interesting cases of non-elementary integrals (see sections \ref{sec:2},  \ref{sec:3} and \ref{sec:4}), we first define the generalized hypergeometric function as it is an important tool that we are going to use in the paper.
\begin{defen}The generalized hypergeometric function, denoted as $_pF_q$, is a special function given by the series \cite{AS,N}
\begin{equation}
_p F_q(a_1, a_2,\cdots,a_p;b_1, b_2, \cdots, b_q; x)=\sum\limits_{n=0}^{\infty}\frac{(a_1)_n (a_2)_n\cdots (a_p)_n}{(b_1)_n (b_2)_n\cdots (b_q)_n}\frac{x^n}{n!},
\label{eq1}
\end{equation} 
where $a_1, a_2,\cdots,a_p$ and $;b_1, b_2, \cdots, b_q$ are arbitrary constants, $(\vartheta)_n=\Gamma(\vartheta+n)/\Gamma(\vartheta)$ (Pochhammer's notation \cite{AS,N}) for any complex $\vartheta$, with $(\vartheta)_0=1$, and $\Gamma$ is the standard gamma function \cite{AS}.
\label{defen2}
\end{defen}

\section{Evaluation of the sine integral $\mbox{Si}_{\beta,\beta+\alpha},\beta\ge1, \alpha>1$}\label{sec:2}
\setcounter{equation}{0}

\begin{teoen}
Let $\beta\ge1$ and $\alpha>1$, and let $\alpha=m\beta+ \epsilon$, where $m$ is an integer ($m\in \mathbb{N}$) and $-\beta<\epsilon<\beta$. 

\begin{enumerate}

\item If $\epsilon=0$, then 
\begin{align}
&\mbox{Si}_{\beta,\beta+\alpha}=\int \frac{\sin{(\lambda x^\beta)}}{\lambda x^{\beta+\alpha}}dx
=\sum\limits_{n=-m}^{-1}(-1)^{n+m}\frac{\lambda^{2n+2m}}{\Gamma(2n+2m+2)}\frac{x^{2\beta n+1}}{2\beta n+1}
\nonumber\\&\hspace{.5cm}+\frac{(-1)^m \lambda^{2m}x}{2^{m+1}\sqrt{\pi}\Gamma(m+1)\Gamma\left(m+\frac{3}{2}\right)(2\beta+1)}\ _2F_3\left(1,1+\frac{1}{2\beta};m+1,m+\frac{3}{2},2+\frac{1}{2\beta};-\frac{\lambda^2 x^{2\beta}}{4}\right)+C, 
\label{eq2}
\end{align}
where $m=\alpha/\beta$.
\item If $\epsilon=1$, then 
\begin{align}
&\mbox{Si}_{\beta,\beta+\alpha}=\int \frac{\sin{(\lambda x^\beta)}}{\lambda x^{\beta+\alpha}}dx
=(-1)^{m}\frac{\lambda^{2m}}{\Gamma(2m+2)}\ln|x|+
\sum\limits_{n=-m}^{-1}(-1)^{n+m}\frac{\lambda^{2n+2m}}{\Gamma(2n+2m+2)}\frac{x^{2\beta n}}{2\beta n}
\nonumber \\&\hspace{.5cm}+\frac{(-1)^{m+1}\lambda^{2m+2}x^{2\beta}}{2^{2m+4}\sqrt{\pi}\Gamma(m+2)\Gamma\left(m+\frac{5}{2}\right)\beta}\ _2F_3\left(1,1;m+2,m+\frac{5}{2},2;-\frac{\lambda^2 x^{2\beta}}{4}\right)+C,
\label{eq3}
\end{align} 
where $m=(\alpha-1)/\beta$.
\item Finally, if $\epsilon \in(-\beta,0)\cup (0,1)\cup(1,\beta)$, we have
\begin{align}
&\int \frac{\sin{(\lambda x^\beta)}}{\lambda x^{\beta+\alpha}}dx
=(-1)^{m}\frac{\lambda^{2m}}{\Gamma(2m+2)}\frac{x^{1-\epsilon}}{1-\epsilon}+
\sum\limits_{n=-m}^{-1}(-1)^{n+m}\frac{\lambda^{2n+2m}}{\Gamma(2n+2m+2)}\frac{x^{2\beta n-\epsilon+1}}{2\beta n-\epsilon+1}
\nonumber \\&\hspace{.5cm}+\frac{(-1)^{m+1}\lambda^{2m+2}x^{2\beta-\epsilon+1}}{2^{2m+3}\sqrt{\pi}\Gamma(m+2)\Gamma\left(m+\frac{5}{2}\right)(2\beta-\epsilon+1)}\ _2F_3\left(1,1+\frac{1-\epsilon}{2\beta};m+2,m+\frac{5}{2},2+\frac{1-\epsilon}{2\beta};-\frac{\lambda^2 x^{2\beta}}{4}\right)+C,
\label{eq4}
\end{align} 
where $m=(\alpha-\epsilon)/\beta$.
\end{enumerate}
\label{teoen:1}
\end{teoen}

\proofen   
We proceed as in \cite{NV,NV1}.
We expand $g(x)$ as Taylor series and integrate the series term by term.
We use the gamma duplication formula\cite{AS}, the gamma property $\Gamma(\alpha+1)=\alpha\Gamma(\alpha)$ 
and Pochhammer's notation (see Definition \ref{defen2}). We also set $\alpha=m\beta+\epsilon$, and then we obtain
\begin{align}
&\int \frac{\sin{(\lambda x^\beta)}}{\lambda x^{\beta+\alpha}}dx
=\int\frac{1}{\lambda x^\beta x^\alpha}\sum\limits_{n=0}^{\infty}(-1)^n\frac{(\lambda x^\beta)^{2n+1}}{(2n+1)!}dx
=\int\sum\limits_{n=0}^{\infty}(-1)^n\frac{\lambda^{2n}}{(2n+1)!}x^{2\beta n-\alpha}dx
\nonumber\\&=\int\sum\limits_{n=0}^{m-1}(-1)^n\frac{\lambda^{2n}}{(2n+1)!}x^{2\beta n-2\beta m-\epsilon}dx
+\int\sum\limits_{n=m}^{\infty}(-1)^n\frac{\lambda^{2n}}{(2n+1)!}x^{2\beta n-2\beta m-\epsilon}dx
\nonumber\\&=\int\sum\limits_{n=0}^{m-1}(-1)^n\frac{\lambda^{2n}}{(2n+1)!}x^{2\beta (n- m)-\epsilon}dx
+\int\sum\limits_{n=m}^{\infty}(-1)^n\frac{\lambda^{2n}}{(2n+1)!}x^{2\beta (n-m)-\epsilon}dx
\nonumber\\&=\int\sum\limits_{n=-m}^{-1}(-1)^{n+m}\frac{\lambda^{2n+2m}}{(2n+2m+1)!}x^{2\beta n-\epsilon}dx
+\int\sum\limits_{n=0}^{\infty}(-1)^{n+m}\frac{\lambda^{2n+2m}}{(2n+2m+1)!}x^{2\beta n-\epsilon}dx
\nonumber\\&=\int\sum\limits_{n=-m}^{-1}(-1)^{n+m}\frac{\lambda^{2n+2m}}{\Gamma(2n+2m+2)}x^{2\beta n-\epsilon}dx
+\int\sum\limits_{n=0}^{\infty}(-1)^{n+m}\frac{\lambda^{2n+2m}}{\Gamma(2n+2m+2)}x^{2\beta n-\epsilon}dx
\label{eq5}\\&=(-1)^{m}\frac{\lambda^{2m}}{\Gamma(2m+2)}\int \frac{dx}{x^{\epsilon}}+
\int\sum\limits_{n=-m}^{-1}(-1)^{n+m}\frac{\lambda^{2n+2m}}{\Gamma(2n+2m+2)}x^{2\beta n-\epsilon}dx 
\nonumber \\&\hspace{.5cm}+ \int\sum\limits_{n=1}^{\infty}(-1)^{n+m}\frac{\lambda^{2n+2m}}{\Gamma(2n+2m+2)}x^{2\beta n-\epsilon}dx 
\nonumber\\&=(-1)^{m}\frac{\lambda^{2m}}{\Gamma(2m+2)}\int \frac{dx}{x^{\epsilon}}+
\int\sum\limits_{n=-m}^{-1}(-1)^{n+m}\frac{\lambda^{2n+2m}}{\Gamma(2n+2m+2)}x^{2\beta n-\epsilon}dx
\nonumber \\&\hspace{.5cm}+ \int\sum\limits_{n=0}^{\infty}(-1)^{n+m+1}\frac{\lambda^{2n+2m+2}}{\Gamma(2n+2m+4)}x^{2\beta n+2\beta-\epsilon}dx
\nonumber\\&=(-1)^{m}\frac{\lambda^{2m}}{\Gamma(2m+2)}\int \frac{dx}{x^{\epsilon}}+
\sum\limits_{n=-m}^{-1}(-1)^{n+m}\frac{\lambda^{2n+2m}}{\Gamma(2n+2m+2)}\frac{x^{2\beta n-\epsilon+1}}{2\beta n-\epsilon+1}
\nonumber \\&\hspace{.5cm}+ \sum\limits_{n=0}^{\infty}(-1)^{n+m+1}\frac{\lambda^{2n+2m+2}}{\Gamma(2n+2m+4)}\frac{x^{2\beta n+2\beta-\epsilon+1}}{2\beta n+2\beta-\epsilon+1}+C_1
\nonumber\\&=(-1)^{m}\frac{\lambda^{2m}}{\Gamma(2m+2)}\int \frac{dx}{x^{\epsilon}}+
\sum\limits_{n=-m}^{-1}(-1)^{n+m}\frac{\lambda^{2n+2m}}{\Gamma(2n+2m+2)}\frac{x^{2\beta n-\epsilon+1}}{2\beta n-\epsilon+1}
\nonumber \\&\hspace{.5cm}+\frac{(-1)^{m+1}\lambda^{2m+2}x^{2\beta-\epsilon+1}}{2^{2m+3}\sqrt{\pi}\Gamma(m+2)\Gamma\left(m+\frac{5}{2}\right)(2\beta-\epsilon+1)}\sum\limits_{n=0}^{\infty}\frac{(1)_n\left(1+\frac{1-\epsilon}{2\beta}\right)_n}{(m+2)_n\left(m+\frac{5}{2}\right)_n\left(2+\frac{1-\epsilon}{2\beta}\right)_n}\frac{\left(-\frac{\lambda^2 x^{2\beta}}{4}\right)^n}{n!}+C_1
\nonumber\\&=(-1)^{m}\frac{\lambda^{2m}}{\Gamma(2m+2)}\int \frac{dx}{x^{\epsilon}}+
\sum\limits_{n=-m}^{-1}(-1)^{n+m}\frac{\lambda^{2n+2m}}{\Gamma(2n+2m+2)}\frac{x^{2\beta n-\epsilon+1}}{2\beta n-\epsilon+1}
\nonumber \\&\hspace{.5cm}+\frac{(-1)^{m+1}\lambda^{2m+2}x^{2\beta-\epsilon+1}}{2^{2m+3}\sqrt{\pi}\Gamma(m+2)\Gamma\left(m+\frac{5}{2}\right)(2\beta-\epsilon+1)}\ _2F_3\left(1,1+\frac{1-\epsilon}{2\beta};m+2,m+\frac{5}{2},2+\frac{1-\epsilon}{2\beta};-\frac{\lambda^2 x^{2\beta}}{4}\right)+C_1
\label{eq6}
\end{align}
\begin{enumerate}
\item For $\epsilon=0$, we substitute $\epsilon=0$ in (\ref{eq5}), and hence, we obtain 
\begin{align}
&\int \frac{\sin{(\lambda x^\beta)}}{\lambda x^{\beta+\alpha}}dx
=\int\frac{1}{\lambda x^\beta x^\alpha}\sum\limits_{n=0}^{\infty}(-1)^n\frac{(\lambda x^\beta)^{2n+1}}{(2n+1)!}dx
\nonumber\\&=\int\sum\limits_{n=-m}^{-1}(-1)^{n+m}\frac{\lambda^{2n+2m}}{\Gamma(2n+2m+2)}x^{2\beta n}dx
+\int\sum\limits_{n=0}^{\infty}(-1)^{n+m}\frac{\lambda^{2n+2m}}{\Gamma(2n+2m+2)}x^{2\beta n}dx
\nonumber\\&=\sum\limits_{n=-m}^{-1}(-1)^{n+m}\frac{\lambda^{2n+2m}}{\Gamma(2n+2m+2)}\frac{x^{2\beta n+1}}{2\beta n+1}
+\sum\limits_{n=0}^{\infty}(-1)^{n+m}\frac{\lambda^{2n+2m}}{\Gamma(2n+2m+2)}\frac{x^{2\beta n+1}}{2\beta n+1}
\nonumber\\&=\sum\limits_{n=-m}^{-1}(-1)^{n+m}\frac{\lambda^{2n+2m}}{\Gamma(2n+2m+2)}\frac{x^{2\beta n+1}}{2\beta n+1}
\nonumber\\&\hspace{.5cm}+\frac{(-1)^m \lambda^{2m}x}{2^{m+1}\sqrt{\pi}(2\beta+1)\Gamma(m+1)\Gamma\left(m+\frac{3}{2}\right)}\sum\limits_{n=0}^{\infty}\frac{(1)_n\left(1+\frac{1}{2\beta}\right)_n}{(m+1)_n\left(m+\frac{3}{2}\right)_n\left(2+\frac{1}{2\beta}\right)_n}\frac{\left(-\frac{\lambda^2 x^{2\beta}}{4}\right)^n}{n!}
\nonumber\\&=\sum\limits_{n=-m}^{-1}(-1)^{n+m}\frac{\lambda^{2n+2m}}{\Gamma(2n+2m+2)}\frac{x^{2\beta n+1}}{2\beta n+1}
\nonumber\\&\hspace{.5cm}+\frac{(-1)^m \lambda^{2m}x}{2^{m+1}\sqrt{\pi}\Gamma(m+1)\Gamma\left(m+\frac{3}{2}\right)(2\beta+1)}\ _2F_3\left(1,1+\frac{1}{2\beta};m+1,m+\frac{3}{2},2+\frac{1}{2\beta};-\frac{\lambda^2 x^{2\beta}}{4}\right)+C
\label{eq7}
\end{align}
which is (\ref{eq2}), and where $m=\alpha/\beta$.
\item For $\epsilon=1$, we set $\epsilon=1$ in (\ref{eq6}) and obtain
\begin{align}
&\int \frac{\sin{(\lambda x^\beta)}}{\lambda x^{\beta+\alpha}}dx
=(-1)^{m}\frac{\lambda^{2m}}{\Gamma(2m+2)}\ln|x|+
\sum\limits_{n=-m}^{-1}(-1)^{n+m}\frac{\lambda^{2n+2m}}{\Gamma(2n+2m+2)}\frac{x^{2\beta n}}{2\beta n}
\nonumber \\&\hspace{.5cm}+\frac{(-1)^{m+1}\lambda^{2m+2}x^{2\beta}}{2^{2m+4}\sqrt{\pi}\Gamma(m+2)\Gamma\left(m+\frac{5}{2}\right)\beta}\ _2F_3\left(1,1;m+2,m+\frac{5}{2},2;-\frac{\lambda^2 x^{2\beta}}{4}\right)+C
\label{eq8}
\end{align} 
which is (\ref{eq3}),  and where $m=(\alpha-1)/\beta$.
\item For $\epsilon\in (-\beta,0)\cup(0,1)\cup(1,\beta)$, (\ref{eq6}) gives
\begin{align}
&\int \frac{\sin{(\lambda x^\beta)}}{\lambda x^{\beta+\alpha}}dx
=(-1)^{m}\frac{\lambda^{2m}}{\Gamma(2m+2)}\frac{x^{1-\epsilon}}{1-\epsilon}+
\sum\limits_{n=-m}^{-1}(-1)^{n+m}\frac{\lambda^{2n+2m}}{\Gamma(2n+2m+2)}\frac{x^{2\beta n-\epsilon+1}}{2\beta n-\epsilon+1}
\nonumber \\&\hspace{.5cm}+\frac{(-1)^{m+1}\lambda^{2m+2}x^{2\beta-\epsilon+1}}{2^{2m+3}\sqrt{\pi}\Gamma(m+2)\Gamma\left(m+\frac{5}{2}\right)(2\beta-\epsilon+1)}\ _2F_3\left(1,1+\frac{1-\epsilon}{2\beta};m+2,m+\frac{5}{2},2+\frac{1-\epsilon}{2\beta};-\frac{\lambda^2 x^{2\beta}}{4}\right)+C
\label{eq9}
\end{align}
which is (\ref{eq4}), and where $m=(\alpha-\epsilon)/\beta$. 
\end{enumerate}
 \hfill$\square$\\[1ex]

\emph{Example~1.} In this example, we evaluate $\int\left[\sin(x^2)/x^{3.5}\right]dx$. We first observe that $\lambda=1$ and $\beta=2$. We also have $3.5=\beta+\alpha=2+1.5=2+(1)2-0.5=\beta+m\beta+\epsilon$, and so $m=1$ and $\epsilon=-0.5$. Substituting $\lambda=1,\beta=2, m=1$ and $\epsilon=-0.5$ in (\ref{eq4}) gives
\begin{equation}
\int\frac{\sin(x^2)}{x^{3.5}}dx=-\frac{x^{1.5}}{9}-\frac{x^{-2.5}}{2.5}+\frac{x^{5.5}}{540\pi}\ _2F_3\left(1,\frac{9}{8};3,\frac{7}{2},\frac{17}{8};-\frac{x^{4}}{4}\right)+C.
\label{eq9-1}
\end{equation}
\\[1ex]

We can use the same procedure for the hyperbolic sine integral, the results are stated in the following theorem. Its proof is similar to that of Theorem \ref{teoen:1}, we will omit it.

\begin{teoen}
Let $\beta\ge1$ and $\alpha>1$, and let $\alpha=m\beta+ \epsilon$, where $m$ is an integer ($m\in \mathbb{N}$) and $-\beta<\epsilon<\beta$. 

\begin{enumerate}

\item If $\epsilon=0$, then 
\begin{align}
&\int \frac{\sinh{(\lambda x^\beta)}}{\lambda x^{\beta+\alpha}}dx
=\sum\limits_{n=-m}^{-1}\frac{\lambda^{2n+2m}}{\Gamma(2n+2m+2)}\frac{x^{2\beta n+1}}{2\beta n+1}
\nonumber\\&\hspace{.5cm}+\frac{ \lambda^{2m}x}{2^{m+1}\sqrt{\pi}\Gamma(m+1)\Gamma\left(m+\frac{3}{2}\right)(2\beta+1)}\ _2F_3\left(1,1+\frac{1}{2\beta};m+1,m+\frac{3}{2},2+\frac{1}{2\beta};\frac{\lambda^2 x^{2\beta}}{4}\right)+C, 
\label{eq10}
\end{align}
where $m=\alpha/\beta$.
\item If $\epsilon=1$, then 
\begin{align}
&\int \frac{\sinh{(\lambda x^\beta)}}{\lambda x^{\beta+\alpha}}dx
=\frac{\lambda^{2m}}{\Gamma(2m+2)}\ln|x|+
\sum\limits_{n=-m}^{-1}\frac{\lambda^{2n+2m}}{\Gamma(2n+2m+2)}\frac{x^{2\beta n}}{2\beta n}
\nonumber \\&\hspace{.5cm}+\frac{\lambda^{2m+2}x^{2\beta}}{2^{2m+4}\sqrt{\pi}\Gamma(m+2)\Gamma\left(m+\frac{5}{2}\right)\beta}\ _2F_3\left(1,1;m+2,m+\frac{5}{2},2;\frac{\lambda^2 x^{2\beta}}{4}\right)+C,
\label{eq11}
\end{align} 
where $m=(\alpha-1)/\beta$.
\item Finally, if $\epsilon \in (-\beta,0)\cup(0,1)\cup(1,\beta)$, we have
\begin{align}
&\int \frac{\sinh{(\lambda x^\beta)}}{\lambda x^{\beta+\alpha}}dx
=\frac{\lambda^{2m}}{\Gamma(2m+2)}\frac{x^{1-\epsilon}}{1-\epsilon}+
\sum\limits_{n=-m}^{-1}\frac{\lambda^{2n+2m}}{\Gamma(2n+2m+2)}\frac{x^{2\beta n-\epsilon+1}}{2\beta n-\epsilon+1}
\nonumber \\&\hspace{.5cm}+\frac{\lambda^{2m+2}x^{2\beta-\epsilon+1}}{2^{2m+3}\sqrt{\pi}\Gamma(m+2)\Gamma\left(m+\frac{5}{2}\right)(2\beta-\epsilon+1)}\ _2F_3\left(1,1+\frac{1-\epsilon}{2\beta};m+2,m+\frac{5}{2},2+\frac{1-\epsilon}{2\beta};\frac{\lambda^2 x^{2\beta}}{4}\right)+C,
\label{eq12}
\end{align} 
where $m=(\alpha-\epsilon)/\beta$.
\end{enumerate}
\label{teoen:2}
\end{teoen}

\section{Evaluation of the cosine integral $\mbox{Ci}_{\beta,2\beta+\alpha},\beta\ge1, \alpha>1$}\label{sec:3}

\begin{teoen}
Let $\beta\ge1$ and $\alpha>1$, and let $\alpha=2\beta m+ \epsilon$, where $m$ is an integer ($m\in \mathbb{N}$) and $ -2\beta<\epsilon<2\beta$. 

\begin{enumerate}

\item If $\epsilon=0$, then 
\begin{align}
&\mbox{Ci}_{\beta,2\beta+\alpha}=\int \frac{\cos{(\lambda x^\beta)}}{\lambda x^{2\beta+\alpha}}dx
=\frac{1}{\lambda}\frac{x^{1-2\beta-\alpha}}{1-2\beta-\alpha}+\sum\limits_{n=-m}^{-1}(-1)^{n+m+1}\frac{\lambda^{2n+2m+1}}{\Gamma(2n+2m+3)}\frac{x^{2\beta n+1}}{2\beta n+1}
\nonumber\\&\hspace{.5cm}+\frac{(-1)^m \lambda^{2m}x}{2^{m+2}\sqrt{\pi}\Gamma\left(m+\frac{3}{2}\right)\Gamma(m+2)(2\beta+1)}\ _2F_3\left(1,1+\frac{1}{2\beta};m+\frac{3}{2},m+2,2+\frac{1}{2\beta};-\frac{\lambda^2 x^{2\beta}}{4}\right)+C,
\label{eq13}
\end{align}
where $m=\alpha/(2\beta)$.
\item If $\epsilon=1$, then 
\begin{align}
&\int \frac{\cos{(\lambda x^\beta)}}{\lambda x^{2\beta+\alpha}}dx
=\frac{1}{\lambda}\frac{x^{1-2\beta-\alpha}}{1-2\beta-\alpha}+\frac{(-1)^{m}\lambda^{2m+1}}{\Gamma(2m+3)}\ln|x|+\sum\limits_{n=-m}^{-1}(-1)^{n+m+1}\frac{\lambda^{2n+2m+1}}{\Gamma(2n+2m+3)}\frac{x^{2\beta n}}{2\beta n}
\nonumber \\&\hspace{3.1cm}+\frac{(-1)^{m+1}\lambda^{2m+3}x^{2\beta}}{2^{2m+5}\sqrt{\pi}\Gamma\left(m+\frac{5}{2}\right)\Gamma(m+3)\beta}\ _2F_3\left(1,1;m+\frac{5}{2},m+3,2;-\frac{\lambda^2 x^{2\beta}}{4}\right)+C,
\label{eq14}
\end{align} 
where $m=(\alpha-1)/(2\beta)$.
\item Finally, if $\epsilon \in  (-2\beta,0)\cup(0,1)\cup(1,2\beta)$, we have
\begin{align}
&\int \frac{\cos{(\lambda x^\beta)}}{\lambda x^{2\beta+\alpha}}dx
=\frac{1}{\lambda}\frac{x^{1-2\beta-\alpha}}{1-2\beta-\alpha}+\frac{(-1)^{m}\lambda^{2m+1}}{\Gamma(2m+3)}\frac{x^{1-\epsilon}}{1-\epsilon}+
\sum\limits_{n=-m}^{-1}(-1)^{n+m+1}\frac{\lambda^{2n+2m+1}}{\Gamma(2n+2m+3)}\frac{x^{2\beta n-\epsilon+1}}{2\beta n-\epsilon+1}
\nonumber \\&\hspace{.5cm}+\frac{(-1)^{m+1}\lambda^{2m+3}x^{2\beta-\epsilon+1}}{2^{2m+4}\sqrt{\pi}\Gamma\left(m+\frac{5}{2}\right)\Gamma(m+3)(2\beta-\epsilon+1)}\ _2F_3\left(1,1+\frac{1-\epsilon}{2\beta};m+\frac{5}{2},m+3,2+\frac{1-\epsilon}{2\beta};-\frac{\lambda^2 x^{2\beta}}{4}\right)+C,
\label{eq15}
\end{align} 
where $m=(\alpha-\epsilon)/(2\beta)$.
\end{enumerate}
\label{teoen:3}
\end{teoen}

\proofen   We proceed as in Theorem \ref{teoen:1}. We have
\begin{align}
&\int \frac{\cos{(\lambda x^\beta)}}{\lambda x^{2\beta+\alpha}}dx
=\int\frac{1}{\lambda x^{2\beta+\alpha} }\sum\limits_{n=0}^{\infty}(-1)^n\frac{(\lambda x^{\beta})^{2n}}{(2n)!}dx
\nonumber\\ &=\int\frac{1}{\lambda x^{2\beta+\alpha}}dx+\frac{1}{\lambda}\int\sum\limits_{n=1}^{\infty}(-1)^n\frac{\lambda^{2n}}{(2n)!}x^{2\beta n-2\beta-\alpha}dx
\nonumber\\ &=\int\frac{1}{\lambda x^{2\beta+\alpha}}dx+\frac{1}{\lambda}\int\sum\limits_{n=0}^{\infty}(-1)^{n+1}\frac{\lambda^{2n+2}}{(2n+2)!}x^{2\beta n-\alpha}dx
\nonumber\\&=\int\frac{1}{\lambda x^{2\beta+\alpha}}dx+\int\sum\limits_{n=0}^{m-1}(-1)^{n+1}\frac{\lambda^{2n+1}}{(2n+2)!}x^{2\beta n-2\beta m-\epsilon}dx
+\int\sum\limits_{n=m}^{\infty}(-1)^{n+1}\frac{\lambda^{2n+1}}{(2n+2)!}x^{2\beta n-2\beta m-\epsilon}dx
\nonumber\\&=\int\frac{1}{\lambda x^{2\beta+\alpha}}dx+\int\sum\limits_{n=0}^{m-1}(-1)^{n+1}\frac{\lambda^{2n+1}}{(2n+2)!}x^{2\beta (n- m)-\epsilon}dx
+\int\sum\limits_{n=m}^{\infty}(-1)^{n+1}\frac{\lambda^{2n+1}}{(2n+2)!}x^{2\beta (n-m)-\epsilon}dx
\nonumber\\&=\int\frac{1}{\lambda x^{2\beta+\alpha}}dx+\int\sum\limits_{n=-m}^{-1}(-1)^{n+m+1}\frac{\lambda^{2n+2m+1}}{(2n+2m+2)!}x^{2\beta n-\epsilon}dx
+\int\sum\limits_{n=0}^{\infty}(-1)^{n+m+1}\frac{\lambda^{2n+2m+1}}{(2n+2m+2)!}x^{2\beta n-\epsilon}dx
\nonumber\\&=\int\frac{1}{\lambda x^{2\beta+\alpha}}dx+\int\sum\limits_{n=-m}^{-1}(-1)^{n+m+1}\frac{\lambda^{2n+2m+1}}{\Gamma(2n+2m+3)}x^{2\beta n-\epsilon}dx
+\int\sum\limits_{n=0}^{\infty}(-1)^{n+m+1}\frac{\lambda^{2n+2m+1}}{\Gamma(2n+2m+3)}x^{2\beta n-\epsilon}dx
\label{eq16}\\&=\int\frac{1}{\lambda x^{2\beta+\alpha}}dx+(-1)^{m+1}\frac{\lambda^{2m+1}}{\Gamma(2m+3)}\int \frac{dx}{x^{\epsilon}}+
\int\sum\limits_{n=-m}^{-1}(-1)^{n+m+1}\frac{\lambda^{2n+2m+1}}{\Gamma(2n+2m+3)}x^{2\beta n-\epsilon}dx 
\nonumber \\&\hspace{.5cm}+ \int\sum\limits_{n=1}^{\infty}(-1)^{n+m+1}\frac{\lambda^{2n+2m+1}}{\Gamma(2n+2m+3)}x^{2\beta n-\epsilon}dx 
\nonumber\\&=\int\frac{1}{\lambda x^{2\beta+\alpha}}dx+(-1)^{m+1}\frac{\lambda^{2m+1}}{\Gamma(2m+3)}\int \frac{dx}{x^{\epsilon}}+
\int\sum\limits_{n=-m}^{-1}(-1)^{n+m+1}\frac{\lambda^{2n+2m+1}}{\Gamma(2n+2m+3)}x^{2\beta n-\epsilon}dx
\nonumber \\&\hspace{.5cm}+ \int\sum\limits_{n=0}^{\infty}(-1)^{n+m}\frac{\lambda^{2n+2m+3}}{\Gamma(2n+2m+5)}x^{2\beta n+2\beta-\epsilon}dx
\nonumber\\&=\frac{1}{\lambda}\frac{x^{1-2\beta-\alpha}}{1-2\beta-\alpha}+(-1)^{m+1}\frac{\lambda^{2m+1}}{\Gamma(2m+3)}\int \frac{dx}{x^{\epsilon}}+
\sum\limits_{n=-m}^{-1}(-1)^{n+m+1}\frac{\lambda^{2n+2m+1}}{\Gamma(2n+2m+3)}\frac{x^{2\beta n-\epsilon+1}}{2\beta n-\epsilon+1}
\nonumber \\&\hspace{.5cm}+\sum\limits_{n=0}^{\infty}(-1)^{n+m}\frac{\lambda^{2n+2m+3}}{\Gamma(2n+2m+5)}\frac{x^{2\beta n+2\beta-\epsilon+1}}{2\beta n+2\beta-\epsilon+1}+C_1
\nonumber\\&=\frac{1}{\lambda}\frac{x^{1-2\beta-\alpha}}{1-2\beta-\alpha}+(-1)^{m+1}\frac{\lambda^{2m+1}}{\Gamma(2m+3)}\int \frac{dx}{x^{\epsilon}}+
\sum\limits_{n=-m}^{-1}(-1)^{n+m+1}\frac{\lambda^{2n+2m+1}}{\Gamma(2n+2m+3)}\frac{x^{2\beta n-\epsilon+1}}{2\beta n-\epsilon+1}
\nonumber \\&\hspace{.5cm}+\frac{(-1)^{m+1}\lambda^{2m+3}x^{2\beta-\epsilon+1}}{2^{2m+4}\sqrt{\pi}\Gamma\left(m+\frac{5}{2}\right)\Gamma(m+3)(2\beta-\epsilon+1)}\sum\limits_{n=0}^{\infty}\frac{(1)_n\left(1+\frac{1-\epsilon}{2\beta}\right)_n}{\left(m+\frac{5}{2}\right)_n(m+3)_n\left(2+\frac{1-\epsilon}{2\beta}\right)_n}\frac{\left(-\frac{\lambda^2 x^{2\beta}}{4}\right)^n}{n!}+C_1
\nonumber\\&=\frac{1}{\lambda}\frac{x^{1-2\beta-\alpha}}{1-2\beta-\alpha}+(-1)^{m}\frac{\lambda^{2m+1}}{\Gamma(2m+3)}\int \frac{dx}{x^{\epsilon}}+
\sum\limits_{n=-m}^{-1}(-1)^{n+m+1}\frac{\lambda^{2n+2m+1}}{\Gamma(2n+2m+3)}\frac{x^{2\beta n-\epsilon+1}}{2\beta n-\epsilon+1}
\nonumber \\&\hspace{.5cm}+\frac{(-1)^{m+1}\lambda^{2m+3}x^{2\beta-\epsilon+1}}{2^{2m+4}\sqrt{\pi}\Gamma\left(m+\frac{5}{2}\right)\Gamma(m+3)(2\beta-\epsilon+1)}\ _2F_3\left(1,1+\frac{1-\epsilon}{2\beta};m+\frac{5}{2},m+3,2+\frac{1-\epsilon}{2\beta};-\frac{\lambda^2 x^{2\beta}}{4}\right)+C_1.
\label{eq17}
\end{align}
\begin{enumerate}
\item For $\epsilon=0$, we substitute $\epsilon=0$ in (\ref{eq16}), and hence, we obtain 
\begin{align}
&\int \frac{\cos{(\lambda x^\beta)}}{\lambda x^{2\beta+\alpha}}dx
\nonumber\\&=\int\frac{dx}{\lambda x^{2\beta+\alpha}}+\int\sum\limits_{n=-m}^{-1}(-1)^{n+m+1}\frac{\lambda^{2n+2m+1}}{\Gamma(2n+2m+3)}x^{2\beta n}dx
\nonumber\\&\hspace{.5cm}+\int\sum\limits_{n=0}^{\infty}(-1)^{n+m+1}\frac{\lambda^{2n+2m+1}}{\Gamma(2n+2m+3)}x^{2\beta n}dx
\nonumber\\&=\frac{1}{\lambda}\frac{x^{1-2\beta-\alpha}}{1-2\beta-\alpha}+\sum\limits_{n=-m}^{-1}(-1)^{n+m+1}\frac{\lambda^{2n+2m+1}}{\Gamma(2n+2m+3)}\frac{x^{2\beta n+1}}{2\beta n+1}
\nonumber\\&\hspace{.5cm}+\sum\limits_{n=0}^{\infty}(-1)^{n+m+1}\frac{\lambda^{2n+2m+1}}{\Gamma(2n+2m+3)}\frac{x^{2\beta n+1}}{2\beta n+1}
\nonumber\\&=\frac{1}{\lambda}\frac{x^{1-2\beta-\alpha}}{1-2\beta-\alpha}+\sum\limits_{n=-m}^{-1}(-1)^{n+m+1}\frac{\lambda^{2n+2m+1}}{\Gamma(2n+2m+3)}\frac{x^{2\beta n+1}}{2\beta n+1}
\nonumber\\&\hspace{.5cm}+\frac{(-1)^m \lambda^{2m}x}{2^{m+2}\sqrt{\pi}(2\beta+1)\Gamma\left(m+\frac{3}{2}\right)\Gamma(m+2)}\sum\limits_{n=0}^{\infty}\frac{(1)_n\left(1+\frac{1}{2\beta}\right)_n}{\left(m+\frac{3}{2}\right)_n(m+2)_n\left(2+\frac{1}{2\beta}\right)_n}\frac{\left(-\frac{\lambda^2 x^{2\beta}}{4}\right)^n}{n!}
\nonumber\\&=\frac{1}{\lambda}\frac{x^{1-2\beta-\alpha}}{1-2\beta-\alpha}+\sum\limits_{n=-m}^{-1}(-1)^{n+m+1}\frac{\lambda^{2n+2m+1}}{\Gamma(2n+2m+3)}\frac{x^{2\beta n+1}}{2\beta n+1}
\nonumber\\&\hspace{.5cm}+\frac{(-1)^m \lambda^{2m}x}{2^{m+2}\sqrt{\pi}\Gamma\left(m+\frac{3}{2}\right)\Gamma(m+2)(2\beta+1)}\ _2F_3\left(1,1+\frac{1}{2\beta};m+\frac{3}{2},m+2,2+\frac{1}{2\beta};-\frac{\lambda^2 x^{2\beta}}{4}\right)+C,
\label{eq18}
\end{align}
which is (\ref{eq13}), and where $m=\alpha/(2\beta)$.
\item For $\epsilon=1$, we set $\epsilon=1$ in (\ref{eq17}) and obtain
\begin{align}
&\int \frac{\cos{(\lambda x^\beta)}}{\lambda x^{2\beta+\alpha}}dx
=\frac{1}{\lambda}\frac{x^{1-2\beta-\alpha}}{1-2\beta-\alpha}+\frac{(-1)^{m}\lambda^{2m+1}}{\Gamma(2m+3)}\ln|x|+\sum\limits_{n=-m}^{-1}(-1)^{n+m+1}\frac{\lambda^{2n+2m+1}}{\Gamma(2n+2m+3)}\frac{x^{2\beta n}}{2\beta n}
\nonumber \\&\hspace{3.1cm}+\frac{(-1)^{m+1}\lambda^{2m+3}x^{2\beta}}{2^{2m+5}\sqrt{\pi}\Gamma\left(m+\frac{5}{2}\right)\Gamma(m+3)\beta}\ _2F_3\left(1,1;m+\frac{5}{2},m+3,2;-\frac{\lambda^2 x^{2\beta}}{4}\right)+C,
\label{eq19}
\end{align} 
which is (\ref{eq14}), and where $m=(\alpha-1)/(2\beta)$.

\item For $\epsilon\in (-2\beta,0)\cup(0,1)\cup(1,2\beta)$, (\ref{eq17}) gives
\begin{align}
&\int \frac{\cos{(\lambda x^\beta)}}{\lambda x^{2\beta+\alpha}}dx
=\frac{1}{\lambda}\frac{x^{1-2\beta-\alpha}}{1-2\beta-\alpha}+\frac{(-1)^{m}\lambda^{2m+1}}{\Gamma(2m+3)}\frac{x^{1-\epsilon}}{1-\epsilon}+
\sum\limits_{n=-m}^{-1}(-1)^{n+m+1}\frac{\lambda^{2n+2m+1}}{\Gamma(2n+2m+3)}\frac{x^{2\beta n-\epsilon+1}}{2\beta n-\epsilon+1}
\nonumber \\&\hspace{.5cm}+\frac{(-1)^{m+1}\lambda^{2m+3}x^{2\beta-\epsilon+1}}{2^{2m+4}\sqrt{\pi}\Gamma\left(m+\frac{5}{2}\right)\Gamma(m+3)(2\beta-\epsilon+1)}\ _2F_3\left(1,1+\frac{1-\epsilon}{2\beta};m+\frac{5}{2},m+3,2+\frac{1-\epsilon}{2\beta};-\frac{\lambda^2 x^{2\beta}}{4}\right)+C,
\label{eq20}
\end{align} 
which is (\ref{eq15}), and where $m=(\alpha-\epsilon)/(2\beta)$.
\end{enumerate}
 \hfill$\square$\\[1ex]

\emph{Example~2.}
In this example, we evaluate $\int\left[\cos(x)/x^{5}\right]dx$. We first observe that $\lambda=1$ and $\beta=1$. We also have $5=2\beta+\alpha=2+3=2+2(1)(1)+1=\beta+2\beta m+\epsilon$, and so $m=1$ and $\epsilon=1$. Substituting $\lambda=1,\beta=1, m=1$ and $\epsilon=1$ in (\ref{eq14}) gives
\begin{equation}
\int\frac{\cos(x)}{x^{5}}dx=-\frac{x^{-4}}{4}-\frac{x^{-2}}{4}+\frac{\ln|x|}{24}+\frac{x^{2}}{720\pi}\ _2F_3\left(1,1;\frac{7}{2},4,2;-\frac{x^{2}}{4}\right)+C.
\label{eq15-1}
\end{equation}
\\[1ex]

We can use the same procedure for the hyperbolic cosine integral, the results are stated in the next theorem. Its proof is similar to Theorem \ref{teoen:3}'s proof, we will omit it.

\begin{teoen}
Let $\beta\ge1$ and $\alpha>1$, and let $\alpha=2\beta m+ \epsilon$, where $m$ is an integer ($m\in \mathbb{N}$) and $ -2\beta<\epsilon<2\beta$. 

\begin{enumerate}

\item If $\epsilon=0$, then 
\begin{align}
&\int \frac{\cosh{(\lambda x^\beta)}}{\lambda x^{2\beta+\alpha}}dx
=\frac{1}{\lambda}\frac{x^{1-2\beta-\alpha}}{1-2\beta-\alpha}+\sum\limits_{n=-m}^{-1}\frac{\lambda^{2n+2m+1}}{\Gamma(2n+2m+3)}\frac{x^{2\beta n+1}}{2\beta n+1}
\nonumber\\&\hspace{.5cm}+\frac{\lambda^{2m}x}{2^{m+2}\sqrt{\pi}\Gamma\left(m+\frac{3}{2}\right)\Gamma(m+2)(2\beta+1)}\ _2F_3\left(1,1+\frac{1}{2\beta};m+\frac{3}{2},m+2,2+\frac{1}{2\beta};\frac{\lambda^2 x^{2\beta}}{4}\right)+C,
\label{eq21}
\end{align}
where $m=\alpha/(2\beta)$.
\item If $\epsilon=1$, then 
\begin{align}
&\int \frac{\cosh{(\lambda x^\beta)}}{\lambda x^{2\beta+\alpha}}dx
=\frac{1}{\lambda}\frac{x^{1-2\beta-\alpha}}{1-2\beta-\alpha}+\frac{\lambda^{2m+1}}{\Gamma(2m+3)}\ln|x|+\sum\limits_{n=-m}^{-1}\frac{\lambda^{2n+2m+1}}{\Gamma(2n+2m+3)}\frac{x^{2\beta n}}{2\beta n}
\nonumber \\&\hspace{3.1cm}+\frac{\lambda^{2m+3}x^{2\beta}}{2^{2m+5}\sqrt{\pi}\Gamma\left(m+\frac{5}{2}\right)\Gamma(m+3)\beta}\ _2F_3\left(1,1;m+\frac{5}{2},m+3,2;\frac{\lambda^2 x^{2\beta}}{4}\right)+C,
\label{eq22}
\end{align} 
where $m=(\alpha-1)/(2\beta)$.
\item Finally, if $\epsilon \in  (-2\beta,0)\cup(0,1)\cup(1,2\beta)$, we have
\begin{align}
&\int \frac{\cos{(\lambda x^\beta)}}{\lambda x^{2\beta+\alpha}}dx
=\frac{1}{\lambda}\frac{x^{1-2\beta-\alpha}}{1-2\beta-\alpha}+\frac{(-1)^{m}\lambda^{2m+1}}{\Gamma(2m+3)}\frac{x^{1-\epsilon}}{1-\epsilon}+
\sum\limits_{n=-m}^{-1}\frac{\lambda^{2n+2m+1}}{\Gamma(2n+2m+3)}\frac{x^{2\beta n-\epsilon+1}}{2\beta n-\epsilon+1}
\nonumber \\&\hspace{.5cm}+\frac{\lambda^{2m+3}x^{2\beta-\epsilon+1}}{2^{2m+4}\sqrt{\pi}\Gamma\left(m+\frac{5}{2}\right)\Gamma(m+3)(2\beta-\epsilon+1)}\ _2F_3\left(1,1+\frac{1-\epsilon}{2\beta};m+\frac{5}{2},m+3,2+\frac{1-\epsilon}{2\beta};\frac{\lambda^2 x^{2\beta}}{4}\right)+C,
\label{eq23}
\end{align} 
where $m=(\alpha-\epsilon)/(2\beta)$.
\end{enumerate}
\label{teoen:4}
\end{teoen}

\section{Evaluation of the exponential integral $\mbox{Ei}_{\beta,\beta+\alpha}, \beta\ge1, \alpha>1$}\label{sec:4}

\begin{teoen}
Let $\beta\ge1$ and $\alpha>1$, and let $\alpha=\beta m+ \epsilon$, where $m$ is an integer ($m\in \mathbb{N}$) and $-\beta<\epsilon<\beta$. 

\begin{enumerate}
\item If $\epsilon=0$, then 
\begin{align}
\mbox{Ei}_{\beta,\beta+\alpha}=\int \frac{e^{\lambda x^\beta}}{\lambda x^{\beta+\alpha}}dx
&=\frac{1}{\lambda}\frac{x^{1-\beta-\alpha}}{1-\beta-\alpha}+\sum\limits_{n=-m}^{-1}\frac{\lambda^{n+m}}{\Gamma(n+m+2)}\frac{x^{\beta n+1}}{\beta n+1}
\nonumber \\&\hspace{0.5cm}+\frac{\lambda^{m}x}{\Gamma(m+2)(\beta+1)}\ _2F_2\left(1,1+\frac{1}{\beta};m+2,2+\frac{1}{\beta};\lambda x^{\beta}\right)+C,
\label{eq24}
\end{align}
where $m=\alpha/(\beta)$.
\item If $\epsilon=1$, then 
\begin{align}
&\mbox{Ei}_{\beta,\beta+\alpha}=\int \frac{e^{\lambda x^\beta}}{\lambda x^{\beta+\alpha}}dx
=\frac{1}{\lambda}\frac{x^{1-\beta-\alpha}}{1-\beta-\alpha}+\frac{\lambda^{m}}{\Gamma(m+2)}\ln|x|+
\sum\limits_{n=-m}^{1}\frac{\lambda^{n+m}}{\Gamma(n+m+2)}\frac{x^{\beta n}}{\beta n}
\nonumber \\&\hspace{4.5cm}+\frac{\lambda^{m+1}x^{\beta}}{\Gamma(m+3)\beta}\ _2F_2\left(1,1;m+3,2;{\lambda x^{\beta}}\right)+C,
\label{eq25}
\end{align} 
where $m=(\alpha-1)/(\beta)$.
\item Finally, if $\epsilon \in  (-\beta,0)\cup(0,1)\cup(1,\beta)$, we have
\begin{align}
&\int \frac{e^{\lambda x^\beta}}{\lambda x^{\beta+\alpha}}dx
=\frac{1}{\lambda}\frac{x^{1-\beta-\alpha}}{1-\beta-\alpha}+\frac{\lambda^{m}}{\Gamma(m+2)}\frac{x^{1-\epsilon}}{1-\epsilon}+
\sum\limits_{n=-m}^{1}\frac{\lambda^{n+m}}{\Gamma(n+m+2)}\frac{x^{\beta n-\epsilon+1}}{\beta n-\epsilon+1}
\nonumber \\&\hspace{2.5cm}+\frac{\lambda^{m+1}x^{\beta-\epsilon+1}}{\Gamma(m+3)(\beta-\epsilon+1)}\ _2F_2\left(1,1+\frac{1-\epsilon}{\beta};m+3,2+\frac{1-\epsilon}{\beta};{\lambda x^{\beta}}\right)+C,
\label{eq26}
\end{align} 
where $m=(\alpha-\epsilon)/(\beta)$.
\end{enumerate}
\label{teoen:5}
\end{teoen}

\proofen   We proceed as before. Then, we have
\begin{align}
&\int \frac{e^{\lambda x^\beta}}{\lambda x^{\beta+\alpha}}dx
=\int\frac{1}{\lambda x^{\beta+\alpha} }\sum\limits_{n=0}^{\infty}\frac{(\lambda x^{\beta})^{n}}{n!}dx
=\int\frac{1}{\lambda x^{\beta+\alpha}}dx+\frac{1}{\lambda}\int\sum\limits_{n=1}^{\infty}\frac{\lambda^{n}}{n!}x^{\beta n-\beta-\alpha}dx
\nonumber\\ &=\int\frac{1}{\lambda x^{\beta+\alpha}}dx+\frac{1}{\lambda}\int\sum\limits_{n=0}^{\infty}\frac{\lambda^{n+1}}{(n+1)!}x^{\beta n-\alpha}dx
\nonumber\\&=\int\frac{1}{\lambda x^{\beta+\alpha}}dx+\int\sum\limits_{n=0}^{m-1}\frac{\lambda^{n}}{(n+1)!}x^{\beta n-\beta m-\epsilon}dx
+\int\sum\limits_{n=m}^{\infty}\frac{\lambda^{n}}{(n+1)!}x^{\beta n-\beta m-\epsilon}dx
\nonumber\\&=\int\frac{1}{\lambda x^{\beta+\alpha}}dx+\int\sum\limits_{n=0}^{m-1}\frac{\lambda^{n}}{(n+1)!}x^{\beta (n- m)-\epsilon}dx
+\int\sum\limits_{n=m}^{\infty}\frac{\lambda^{n}}{(n+1)!}x^{\beta (n-m)-\epsilon}dx
\nonumber\\&=\int\frac{1}{\lambda x^{\beta+\alpha}}dx+\int\sum\limits_{n=-m}^{-1}\frac{\lambda^{n+m}}{(n+m+1)!}x^{\beta n-\epsilon}dx
+\int\sum\limits_{n=0}^{\infty}\frac{\lambda^{n+m}}{(n+m+1)!}x^{\beta n-\epsilon}dx
\nonumber\\&=\int\frac{1}{\lambda x^{\beta+\alpha}}dx+\int\sum\limits_{n=-m}^{-1} \frac{\lambda^{n+m}}{\Gamma(n+m+2)}x^{\beta n-\epsilon}dx
+\int\sum\limits_{n=0}^{\infty}\frac{\lambda^{n+m}}{\Gamma(n+m+2)}x^{\beta n-\epsilon}dx
\label{eq27}\\&=\int\frac{1}{\lambda x^{\beta+\alpha}}dx+\frac{\lambda^{m}}{\Gamma(m+2)}\int \frac{dx}{x^{\epsilon}}+
\int\sum\limits_{n=-m}^{-1}\frac{\lambda^{n+m}}{\Gamma(n+m+2)}x^{\beta n-\epsilon}dx 
+ \int\sum\limits_{n=1}^{\infty}\frac{\lambda^{n+m}}{\Gamma(n+m+2)}x^{\beta n-\epsilon}dx 
\nonumber\\&=\int\frac{1}{\lambda x^{\beta+\alpha}}dx+\frac{\lambda^{m}}{\Gamma(m+2)}\int \frac{dx}{x^{\epsilon}}+
\int\sum\limits_{n=-m}^{-1}\frac{\lambda^{n+m}}{\Gamma(n+m+2)}x^{\beta n-\epsilon}dx
+ \int\sum\limits_{n=0}^{\infty}\frac{\lambda^{n+m+1}}{\Gamma(n+m+3)}x^{\beta n+\beta-\epsilon}dx
\nonumber\\&=\frac{1}{\lambda}\frac{x^{1-\beta-\alpha}}{1-\beta-\alpha}+\frac{\lambda^{m}}{\Gamma(m+2)}\int \frac{dx}{x^{\epsilon}}+
\sum\limits_{n=-m}^{-1}\frac{\lambda^{n+m}}{\Gamma(n+m+2)}\frac{x^{\beta n-\epsilon+1}}{\beta n-\epsilon+1}
+\sum\limits_{n=0}^{\infty}\frac{\lambda^{n+m+1}}{\Gamma(n+m+3)}\frac{x^{\beta n+\beta-\epsilon+1}}{\beta n+\beta-\epsilon+1}+C_1
\nonumber\\&=\frac{1}{\lambda}\frac{x^{1-\beta-\alpha}}{1-\beta-\alpha}+\frac{\lambda^{m}}{\Gamma(m+2)}\int \frac{dx}{x^{\epsilon}}+
\sum\limits_{n=-m}^{-1}\frac{\lambda^{n+m}}{\Gamma(n+m+2)}\frac{x^{\beta n-\epsilon+1}}{\beta n-\epsilon+1}
\nonumber \\&\hspace{5.65cm}+\frac{\lambda^{m+1}x^{\beta-\epsilon+1}}{\Gamma(m+3)(\beta-\epsilon+1)}\sum\limits_{n=0}^{\infty}\frac{(1)_n\left(1+\frac{1-\epsilon}{\beta}\right)_n}{(m+3)_n\left(2+\frac{1-\epsilon}{\beta}\right)_n}\frac{\left({\lambda x^{\beta}}\right)^n}{n!}+C_1
\nonumber\\&=\frac{1}{\lambda}\frac{x^{1-\beta-\alpha}}{1-\beta-\alpha}+\frac{\lambda^{m}}{\Gamma(m+2)}\int \frac{dx}{x^{\epsilon}}+
\sum\limits_{n=-m}^{-1}\frac{\lambda^{n+m}}{\Gamma(n+m+2)}\frac{x^{\beta n-\epsilon+1}}{\beta n-\epsilon+1}
\nonumber \\&\hspace{4.5cm}+\frac{\lambda^{m+1}x^{\beta-\epsilon+1}}{\Gamma(m+3)(\beta-\epsilon+1)}\ _2F_2\left(1,1+\frac{1-\epsilon}{\beta};m+3,2+\frac{1-\epsilon}{\beta};{\lambda x^{\beta}}\right)+C_1.
\label{eq28}
\end{align}
\begin{enumerate}
\item For $\epsilon=0$, we substitute $\epsilon=0$ in (\ref{eq27}), and hence, we obtain 
\begin{align}
&\int \frac{e^{\lambda x^\beta}}{\lambda x^{\beta+\alpha}}dx
=\int\frac{dx}{\lambda x^{\beta+\alpha}}+\int\sum\limits_{n=-m}^{-1}\frac{\lambda^{n+m}}{\Gamma(n+m+2)}x^{\beta n}dx
+\int\sum\limits_{n=0}^{\infty}\frac{\lambda^{n+m}}{\Gamma(n+m+2)}x^{\beta n}dx
\nonumber\\&=\frac{1}{\lambda}\frac{x^{1-\beta-\alpha}}{1-\beta-\alpha}+\sum\limits_{n=-m}^{-1}\frac{\lambda^{n+m}}{\Gamma(n+m+2)}\frac{x^{\beta n+1}}{\beta n+1}
+\sum\limits_{n=0}^{\infty}\frac{\lambda^{n+m}}{\Gamma(n+m+2)}\frac{x^{\beta n+1}}{\beta n+1}
\nonumber\\&=\frac{1}{\lambda}\frac{x^{1-\beta-\alpha}}{1-\beta-\alpha}+\sum\limits_{n=-m}^{-1}\frac{\lambda^{n+m}}{\Gamma(n+m+2)}\frac{x^{\beta n+1}}{\beta n+1}
+\frac{\lambda^{m}x}{\Gamma(m+2)(\beta+1)}\sum\limits_{n=0}^{\infty}\frac{(1)_n\left(1+\frac{1}{\beta}\right)_n}{(m+2)_n\left(2+\frac{1}{\beta}\right)_n}\frac{\left({\lambda x^{\beta}}\right)^n}{n!}
\nonumber\\&=\frac{1}{\lambda}\frac{x^{1-\beta-\alpha}}{1-\beta-\alpha}+\sum\limits_{n=-m}^{-1}\frac{\lambda^{n+m}}{\Gamma(n+m+2)}\frac{x^{\beta n+1}}{\beta n+1}
+\frac{\lambda^{m}x}{\Gamma(m+2)(\beta+1)}\ _2F_2\left(1,1+\frac{1}{\beta};m+2,2+\frac{1}{\beta};\lambda x^{\beta}\right)+C,
\label{eq29}
\end{align}
which is (\ref{eq24}), and where $ m=\alpha/\beta$.

\item For $\epsilon=1$, we set $\epsilon=1$ in (\ref{eq28}) and obtain
\begin{align}
&\int \frac{e^{\lambda x^\beta}}{\lambda x^{\beta+\alpha}}dx
=\frac{1}{\lambda}\frac{x^{1-\beta-\alpha}}{1-\beta-\alpha}+\frac{\lambda^{m}}{\Gamma(m+2)}\ln|x|+
\sum\limits_{n=-m}^{-1}\frac{\lambda^{n+m}}{\Gamma(n+m+2)}\frac{x^{\beta n}}{\beta n}
\nonumber \\&\hspace{2.5cm}+\frac{\lambda^{m+1}x^{\beta}}{\Gamma(m+3)\beta}\ _2F_2\left(1,1;m+3,2;{\lambda x^{\beta}}\right)+C,
\end{align} 
which is (\ref{eq25}), and where $ m=(\alpha-1)/\beta$.
\item For $\epsilon\in (-\beta,0)\cup(0,1)\cup(1,\beta)$, (\ref{eq28}) gives
\begin{align}
&\int \frac{e^{\lambda x^\beta}}{\lambda x^{\beta+\alpha}}dx
=\frac{1}{\lambda}\frac{x^{1-\beta-\alpha}}{1-\beta-\alpha}+\frac{\lambda^{m}}{\Gamma(m+2)}\frac{x^{1-\epsilon}}{1-\epsilon}+
\sum\limits_{n=-m}^{-1}\frac{\lambda^{n+m}}{\Gamma(n+m+2)}\frac{x^{\beta n-\epsilon+1}}{\beta n-\epsilon+1}
\nonumber \\&\hspace{2.5cm}+\frac{\lambda^{m+1}x^{\beta-\epsilon+1}}{\Gamma(m+3)(\beta-\epsilon+1)}\ _2F_2\left(1,1+\frac{1-\epsilon}{\beta};m+3,2+\frac{1-\epsilon}{\beta};{\lambda x^{\beta}}\right)+C,
\label{eq30}
\end{align} 
which is (\ref{eq26}), and where $ m=(\alpha-\epsilon)/\beta$.
\end{enumerate}
 \hfill$\square$\\[1ex]

\emph{Example~3.}
In this example, we evaluate $\int(e^{-x^2}/x^{4})dx$. We first observe that $\lambda=-1$ and $\beta=2$. We also have $4=\beta+\alpha=2+2=2+2(1)+0=\beta+\beta m+\epsilon$, and so $m=1$ and $\epsilon=0$. Substituting $\lambda=1,\beta=1, m=1$ and $\epsilon=0$ in (\ref{eq24}) gives
\begin{equation}
\int\frac{e^{-x^2}}{x^{4}}dx=\frac{x^{-3}}{3}-\frac{1}{x}-\frac{x}{4}\ _2F_2\left(1,2;3,3;-x^{2}\right)+C.
\label{eq30-1}
\end{equation}
\\[1ex]

\begin{coren}
Let  $\alpha>1$ and let $\alpha=m+ \epsilon$, where $m$ is an integer ($m\in \mathbb{N}$) and $-1<\epsilon\le1$. 

\begin{enumerate}
\item If $\epsilon=0$ or $1$, then 
\begin{align}
&\mbox{Ei}_{1,1+\alpha}=\int \frac{e^{\lambda x}}{\lambda x^{1+\alpha}}dx
=-\frac{1}{\lambda \alpha x^{\alpha}}+\frac{\lambda^{m}}{\Gamma(m+2)}\ln|x|+
\sum\limits_{n=-m}^{-1}\frac{\lambda^{n+m}}{\Gamma(n+m+2)}\frac{x^{n}}{n}
\nonumber \\&\hspace{4.5cm}+\frac{\lambda^{m+1}x}{\Gamma(m+3)\beta}\ _2F_2\left(1,1;m+3,2;{\lambda x}\right)+C,
\label{eq31}
\end{align} 
where $m=\alpha-1$.
\item And if $\epsilon \in (-1,0)\cup(0,1)$, we have
\begin{align}
&\int \frac{e^{\lambda x}}{\lambda x^{1+\alpha}}dx
=-\frac{1}{\lambda\alpha x^{\alpha}}+\frac{\lambda^{m}}{\Gamma(m+2)} \frac{x^{1-\epsilon}}{1-\epsilon}+
\sum\limits_{n=-m}^{-1}\frac{\lambda^{n+m}}{\Gamma(n+m+2)}\frac{x^{n-\epsilon+1}}{n-\epsilon+1}
\nonumber \\&\hspace{2.5cm}+\frac{\lambda^{m+1}x^{2-\epsilon}}{\Gamma(m+3)(2-\epsilon)}\ _2F_2\left(1,2-\epsilon;m+3,3-\epsilon;{\lambda x}\right)+C,
\label{eq32}
\end{align} 
where $m=\alpha-\epsilon$.
\end{enumerate}
\label{coren:1}
\end{coren}

\proofen
\begin{enumerate}
\item If $\epsilon=0$ or $1$ implies $\alpha=m+\epsilon$ is an integer ($\alpha\in \mathbb{N}$) since  ($m\in \mathbb{N}$). Morever,  $\alpha=m+\epsilon$ implies $\beta=1$ in Theorem \ref{teoen:5}. Therefore, we obtain (\ref{eq31}) by setting $\beta=1$ in (\ref{eq25}).
\item For $\epsilon \in  (-1,0)\cup(0,1)$, we set $\beta=1$ in (\ref{eq26}) and obtain (\ref{eq32}).
\end{enumerate}
\hfill $\square$

\emph{Example~4.}
In this example, we evaluate $\int\left(e^{-x}/x^{3.7}\right)dx$. We first observe that $\lambda=-1$. We also have $3.7=1+\alpha=1+2.7=1+2+0.7=1+m+\epsilon$, and so $m=2$ and $\epsilon=0.7$. Substituting $\lambda=-1, m=2$ and $\epsilon=0.7$ in (\ref{eq32}) gives
\begin{equation}
\int\frac{e^{-x}}{x^{3.7}}dx=\frac{x^{-2.7}}{2.7}-\frac{x^{0.3}}{1.8}-\frac{x^{-1.7}}{1.7}+\frac{x^{-0.7}}{1.4}-\frac{x^{1.3}}{31.2}\ _2F_2\left(1,1.3;5,2.3;-x\right)+C.
\label{eq32-1}
\end{equation}
\\[1ex]

\section{Conclusion}\label{sec:5}\setcounter{equation}{0}
Formulas for the non-elementary integrals $\mbox{Si}_{\beta,\alpha}=\int [\sin{(\lambda x^\beta)}/{(\lambda x^\alpha)}] dx, \beta\ge1, \alpha>\beta+1$, and $\mbox{Ci}_{\beta,\alpha}=\int [\cos{(\lambda x^\beta)}/{(\lambda x^\alpha)}] dx,\beta\ge1, \alpha>2\beta+1$, were explicitly derived in terms of the hypergeometric function ${}_2F_3$ (see Theorems \ref{teoen:1} and \ref{teoen:2}). Once derived, formulas for the hyperbolic sine and hyperbolic cosine integrals were deduced from those of the sine and cosine integrals  (see Theorems \ref{teoen:2} and \ref{teoen:4}).
On the other hand, the exponential integral $\mbox{Ei}_{\beta,\alpha}=\int (e^{\lambda x^\beta}/x^\alpha)dx, \beta\ge1, \alpha> \beta+1$  was expressed in terms of the hypergeometric function ${}_2F_2$ (see Theorem \ref{teoen:5} and Corollary  \ref{coren:1}). 

Beside, illustrative examples were given. Therefore, their corresponding definite integrals can now be evaluated using the FTC rather than using numerical integration.

\begin{Biblioen}
\bibitem{AS}{\bf Abramowitz~M., Stegun~I.A.} Handbook of mathematical functions with formulas,
graphs and mathematical tables. National Bureau of Standards,1964. 1046~p.
\bibitem{C}{\bf Chiccoli~C., Lorenzutta~S., Maino~G.} Concerning some integrals of the generalized exponential-integral function//Computers Math. Applic., 1992. Vol.~23, No.~11, P.~13--21.
\bibitem{Ch}{\bf Chen~X.} Exponential asymptotics and law of the iterated logarithm for intersection local times of random walks// Ann. Probab., 2004. Vol.~32, No.~4, P.~3248--3300. DOI 10.1214/009117904000000513
\bibitem{MZ}{\bf Marchisotto~E.A., Zakeri~G.-A.} An invitation to integration in finite terms// College
Math. J., 1994. Vol.~25, no~4. P.~295--308. DOI:~10.2307/2687614
\bibitem{NV}{\bf Nijimbere~V.} Evaluation of the non-elementary integral $\int e^{\lambda x^\alpha} dx$, $\alpha\ge2$, and other related integrals// Ural Math. J., 2017. Vol~3, no.~2. P.~130--142. DOI:~10.15826/Umj.2017.2.014
\bibitem{NV1}{\bf Nijimbere~V.} Evaluation of some non-elementary integrals involving sine, cosine, exponential and logarithmic integrals: Part I// Ural Math. J., 2017. Accepted for publication. 
\bibitem{N} NIST Digital Library of Mathematical Functions. \url{http://dlmf.nist.gov/}
\bibitem{Ra}{\bf Rahman~M.} Applications of Fourier transforms to generalized functions. Witt Press, 2011. 192~p.
\bibitem{R}{\bf Rosenlicht~M.} Integration in finite terms// Amer. Math. Monthly, 1972. Vol~79, no.~9. P.~963--972. 
DOI:~10.2307/2318066
\bibitem{S}{\bf Shore~S.N.} Blue Sky and Hot Piles: The Evolution of Radiative Transfer Theory from Atmospheres to Nuclear Reactors// Historia Mathematica, 2002. Vol~29, P.~463--489. Doi:10.1006/hmat.2002.2360
\end{Biblioen}

 \end{document}